\documentclass[11pt]{article}
\usepackage{listings}

\usepackage{pdflscape}
\usepackage{epsfig}
\usepackage{lscape}
\usepackage{longtable}
\usepackage{amsmath,amssymb}
\usepackage{graphicx}
\usepackage{color}
\usepackage{placeins}
\usepackage{url}
\def\sqr#1#2{{\vcenter{\vbox{\hrule height.#2pt
        \hbox{\vrule width.#2pt height#1pt \kern#2pt
        \vrule width.#2pt}
        \hrule height.#2pt}}}}

\def\approxleq{ \kern3pt \mbox{\raisebox{.6ex}{$<$}} \kern-8pt
  \mbox{\raisebox{-.6ex}{$\sim$}} \kern5pt}

\def\norm#1{\|#1 \|}
\def\inprod#1#2{\langle#1,\,#2 \rangle}

\def\cA{{\cal A}}   \def\cK{{\cal K}}
\def\cS{{\cal S}}
\def\cX{{\cal X}} \def\cZ{{\cal Z}}
\def\cQ{{\cal Q}}

\def\nn{\nonumber}

\def\Xprime{X^\prime}

\def\cL{{\cal L}}

\def \[{\begin{equation}}
\def \]{\end{equation}}

\newlength{\len}
\newtheorem{theorem}{Theorem}[section]

\newtheorem{lemma}{Lemma}[section]

\newtheorem{remark}{Remark}[section]
\newtheorem{assumption}{Assumption}[section]

\DeclareMathOperator*{\argmin}{argmin}

\begin{document}
\title{A Convergent  $3$-Block Semi-Proximal ADMM for Convex Minimization
Problems with One Strongly Convex Block}

\author{Min Li\thanks{School of Economics and Management, Southeast
University, Nanjing, 210096, China ({\tt limin@seu.edu.cn}). This
author was supported by the National Natural Science Foundation of
China (Grant No. 11001053), Program for New Century Excellent
Talents in University (Grant No. NCET-12-0111) and Qing Lan
Project.}, \; Defeng Sun\thanks{Department  of  Mathematics  and
Risk  Management Institute, National University of Singapore, 10
Lower Kent Ridge Road, Singapore ({\tt matsundf@nus.edu.sg}). } \
and \ Kim-Chuan Toh\thanks{Department of Mathematics, National
University of Singapore, 10 Lower Kent Ridge Road, Singapore ({\tt
mattohkc@nus.edu.sg}).  } }
\date{October 29, 2014}

\maketitle

\begin{abstract}
In this paper, we present a  semi-proximal alternating direction
method of multipliers (ADMM) for solving $3$-block separable convex
minimization problems  with the second block  in the objective being
a strongly convex  function and one  coupled linear equation
constraint.  By choosing the semi-proximal terms properly, we
establish the global convergence of the proposed semi-proximal ADMM
for the step-length $\tau \in (0, (1+\sqrt{5})/2)$ and the  penalty
parameter $\sigma\in (0, +\infty)$. In particular, if $\sigma>0$ is
smaller than a certain threshold and the first and third  linear
operators in the linear equation constraint are injective, then all
the three added semi-proximal terms can be dropped and consequently,
the convergent $3$-block semi-proximal ADMM reduces to the directly
extended $3$-block ADMM with  $\tau \in (0, (1+\sqrt{5})/2)$.
\end{abstract}

\medskip
{\small
\begin{center}
\parbox{0.95\hsize}{{\bf Keywords.}\;
Convex minimization
problems, alternating direction method of multipliers, semi-proximal, strongly convex.}
\end{center}

\begin{center}
\parbox{0.95\hsize}{{\bf AMS subject classifications.}\; 90C25, 90C33, 65K05}
\end{center}}

\section{Introduction}\label{Introduction}

We consider the following separable convex minimization problem  whose  objective function is the sum of three
functions without coupled variables:
\[  \label{ConvexP-G}
    \min_{x_1, x_2, x_3} \Bigl\{  \theta_1(x_1)  + \theta_2(x_2)  +\theta_3(x_3) \; \Big| \;
    A_1^* x_1 +  A_2^* x_2+  A_3^* x_3 = c, \; x_i \in {\cal X}_i, \; i=1,2,3  \Bigr\},
   \]
where ${\cal X}_i$ ($i = 1, 2, 3$) and ${\cal Z}$ are real finite
dimensional Euclidean spaces each equipped with an inner product
$\langle \cdot, \cdot\rangle$ and its induced norm $\|\cdot\|$,
 $\theta_i: {\cal X}_i \rightarrow (-\infty, +\infty]$ ($i = 1,
2, 3$) are closed proper convex functions, $A_i^*: {\cal X}_i
\rightarrow {\cal Z}$   is the adjoint of the linear operator
$A_i:\cZ \to \cX_i$, $i=1,2,3$,  and $c \in {\cal Z}$. Since
$\theta_i$,  $i=1, 2, 3$, are closed proper convex functions, there
exist self-adjoint and positive semi-definite operators  $\Sigma_i$,
  $i = 1, 2, 3$, such that
\[\label{strong-convex}
       \big\langle \hat{x}_i - x_i, \; \hat{w}_i - w_i\big\rangle \ge \langle \hat{x}_i - x_i, {\Sigma_i}(\hat{x}_i - x_i)\rangle
       \quad \forall\; \hat{x}_i, x_i \in  {\rm dom}(\theta_i), \; \hat{w}_i \in \partial \theta_i(\hat{x}_i),
       w_i \in \partial \theta_i(x_i), \]
where $\partial  \theta_i$ is the sub-differential mapping  of
$\theta_i$, $i=1, 2, 3$. The solution set of problem
\eqref{ConvexP-G} is assumed to be nonempty throughout our
discussions in this paper.

Let $\sigma > 0$ be a given penalty parameter and $z \in {\cal Z}$
be the Lagrange multiplier associated with the linear equality
constraint in  problem \eqref{ConvexP-G}. For any $(x_1, x_2, x_3) \in  {\cal X}_1\times {\cal X}_2\times {\cal X}_3$, write
$x\equiv(x_1, x_2, x_3) $, $\theta (x)\equiv  \theta_1(x_1)  + \theta_2(x_2)  +\theta_3(x_3) $ and $A^*x \equiv A_1^* x_1 +  A_2^* x_2+  A_3^* x_3$.
Then the augmented Lagrangian function
for problem  \eqref{ConvexP-G} is defined by
\[\label{Lagrange-func2} {\cal L}_{\sigma}(x_1,x_2, x_3;z) : =
   \theta (x) + \langle z, \;   A^* x  - c\rangle +
 \frac{\sigma}{2}\|  A^* x  - c\|^2  \]
for any  $(x_1,x_2, x_3, z) \in {\cal X}_1 \times {\cal X}_2 \times
{\cal X}_3 \times {\cal Z}$. The direct extension of the classical
alternating direction method of multipliers (ADMM) for solving
problem \eqref{ConvexP-G} consists of the following iterations for
$k=0,1, \ldots$
\begin{equation}\label{extend-ADMM}
\left\{
\begin{array}{l}
\displaystyle  x_1^{k+1} := \argmin_{x_1 \in {\cal X}_1}  \bigl\{
     {\cal L}_{\sigma}(x_1, x_2^k, x_3^k; z^k)   \bigr\}, \\
\displaystyle x_2^{k+1} := \argmin_{x_2 \in {\cal X}_2}  \bigl\{
     {\cal L}_{\sigma}(x_1^{k+1}, x_2, x_3^k; z^k) \bigr\}, \\
\displaystyle x_3^{k+1} := \argmin_{x_3 \in {\cal X}_3}  \bigl\{
     {\cal L}_{\sigma}(x_1^{k+1}, x_2^{k+1},  x_3; z^k)
        \bigr\}, \\
z^{k+1} := z^k + \tau\sigma(  A^* x^{k+1} - c),
\end{array}
\right.
\end{equation}
where $\tau >0$ is the step-length.
Different from the $2$-block ADMM whose convergence has been established for a long time
\cite{glowinski1975sur, gabay1976dual, Glowinski1980lectures, fortin1983augmentedlag,Gabay1983299, eckstein1992douglas},
the $3$-block ADMM  may not converge in general, which  was
demonstrated   by  Chen, He, Ye and Yuan \cite{ChenHeYeYuan} using
counterexamples.  Nevertheless, if all the functions $\theta_i$,
$i=1,2,3$, are strongly convex,
 Han and Yuan \cite{HanYuan2012} proved the
global convergence of  the $3$-block ADMM  scheme \eqref{extend-ADMM} with
$\tau=1$ (Han and Yuan actually considered the general  $m$-block case for
any $m\ge 3$. Here and below we focus on the $3$-block case only) under the
condition that
\[  \Sigma_i = \mu_i I \succ 0, \; i = 1, 2,3, \quad 0 < \sigma \le
  \min_{i=1, 2, 3}\Big\{\frac{\mu_i}{3\lambda_{\max}(A_iA_i^*)}\Big\}, \nn \]
  where $\lambda_{\max}(S)$ is the largest eigenvalue of a given self-adjoint linear operator $ S $.
Hong and Luo \cite{HongLuo2012} proposed to adopt a small
step-length $\tau$ when updating the Lagrange multiplier $z^{k+1}$
in \eqref{extend-ADMM}. Chen, Shen and You \cite{ChenShenYou}
proposed the following sufficient condition
\[  A_1^* \; \hbox{is injective}, \;\; \Sigma_i = \mu_i I \succ 0,
\; i = 2, 3 \;\;  \hbox{and} \;\; 0 < \sigma \le  \min\Big\{\frac{\mu_2}{\lambda_{\max}(A_2A_2^*)},
\; \frac{\mu_3}{\lambda_{\max}(A_3A_3^*)} \Big\} \nn \]
 for the global convergence of the directly extended   $3$-block ADMM
 with $\tau = 1$ for solving problem \eqref{ConvexP-G}.  Closely related to the work of  Chen, Shen and You \cite{ChenShenYou},
in \cite{LinMaZhang2014A}, Lin, Ma and Zhang provided an analysis on the
 iteration complexity  for the same method   under the
condition
\[ \Sigma_i = \mu_i I \succ 0, \; i = 2, 3 \;\;  \hbox{and} \;\; 0 < \sigma \le
\min\Big\{\frac{\mu_2}{2\lambda_{\max}(A_2A_2^*)}, \;
\frac{\mu_3}{2\lambda_{\max}(A_3A_3^*)} \Big\}.  \nn \]
 In \cite{LinMaZhang2014B}, under additional assumptions including some smoothness conditions, the same group of
 authors further  proved the global linear convergence of the  mentioned method.

The purpose of this work is to extend  the $2$-block semi-proximal
ADMM studied in \cite{FsPoSunTse} to deal with problem
\eqref{ConvexP-G}  by  only  assuming   $\theta_2$ to be strongly
convex, i.e., $\Sigma_2\succ 0$. Note that the semi-proximal ADMM
with $\tau >1$ often works better
 in practice than its counterpart with  $\tau \le 1$.  So it is desirable  to establish the convergence of the proposed semi-proximal
ADMM that allows  $\tau$ to stay in the larger region
$ (0, (1 + \sqrt{5})/2)$.

One of our motivating examples is  the following convex quadratic conic
programming
\begin{eqnarray}
  \begin{array}{ll}
    \min &  \displaystyle \frac{1}{2} \inprod{X}{\cQ X} + \inprod{C}{X}  \\[5pt]
   \mbox{s.t.}
       &  \cA  X    \ge  b, \quad
       X \in  \cK  \, ,
\end{array}
 \label{eq-qsdp}
\end{eqnarray}
where  $\cK$ is a nonempty  closed convex cone in  a finite dimensional real Euclidean space $\cX$
  endowed with an  inner
product $\inprod{\cdot}{\cdot}$ and its induced norm
$\norm{\cdot}$, $\cQ:\cX \to \cX$ is a self-adjoint and positive semi-definite linear
operator,
 $\cA:\cX \rightarrow \Re^{m}$ is a linear map,  $C\in \cX$ and $b \in \Re^{m}$
 are given data. The dual of problem  (\ref{eq-qsdp})  takes the form of
 \begin{equation}   \label{eq-d-qsdp}
  \begin{array}{ll}
    \max &   \displaystyle  - \frac{1}{2}\inprod{\Xprime}{\cQ \Xprime}   + \inprod{b}{y}   \\[5pt]
   \mbox{s.t.} &  \cA^*y - \cQ \Xprime  + S = C, \quad
                y \geq 0, \quad S \in \cK^* \, ,
   \end{array}
\end{equation}
where $\cK^*:=\{v \in \cX: \langle v, w \rangle \ge 0 \; \forall\,  w
\in \cK\}$ is the dual cone of $\cK$.
Since $\cQ$ is  self-adjoint and positive semi-definite, $\cQ$ can be decomposed as
$\cQ = \cL^*\cL$ for some linear map $\cL$. By introducing a new
variable $\Xi = -\cL\Xprime$, we can re-write   problem \eqref{eq-d-qsdp} equivalently
as
 \begin{eqnarray}
  \begin{array}{rllll}
    \min & \displaystyle  \delta_{\Re^{m}_+}(y)- \inprod{b}{y}   + \frac{1}{2}\norm{\Xi}^2 + \delta_{\cK^*}(S)      \\[5pt]
   \mbox{s.t.} &   \cA^*y + \cL^*\Xi +  S   = C,
   \end{array}
   \label{eq-d1-qsdp-eg}
\end{eqnarray}
where $\delta_{\Re^{m}_+}(\cdot)$ and $\delta_{\cK^*}(\cdot)$ are the indicator functions of $\Re^{m}_+$ and $\cK^*$, respectively.  As one can see, problem \eqref{eq-d1-qsdp-eg} has only one strongly convex block, i.e., the block with respect to $\Xi$. Consequently, the results in the aforementioned papers for the convergence analysis of   the directly extended   3-block ADMM   applied to solving problem \eqref{eq-d1-qsdp-eg} are no longer valid. We shall show in the next section that our proposed 3-block semi-proximal ADMM
 can exactly  solve this kind of problems.
When $\cK = \cS_+^n$, the cone of symmetric and positive
semi-definite matrices in the space $\cS^n$ of $n\times n$ symmetric
matrices, problem \eqref{eq-d1-qsdp-eg} is a convex quadratic
semidefinite programming problem that has been extensively studied
both theoretically and numerically in the literature \cite{LiSunToh,
NieYuanAOR,Qi09MOR,QiSunIMA11,SunMOR06,SunSunZhangMP08,SunZhang10EJOR, TohMP08,
TTTPJO07, ZhaoXY09}, to name only a few.

The remaining parts of this paper are organized as follows.
In the next section, we first  present our $3$-block semi-proximal ADMM and then
provide the main convergence results. We  give some concluding remarks in
the final section.

\medskip

{\bf Notation.} \begin{itemize}
\item The effective domain of a function $f$: ${\cal X}
\rightarrow (-\infty, +\infty]$ is defined as $\hbox{dom}(f): = \{x
\in {\cal X} \; | \; f(x) < + \infty\}$. The set of all relative
interior points of a given nonempty convex set ${\cal C}$    is denoted by
{\rm ri}$(\cal C)$.

\item  For convenience,  for any given $x$, we use $\|x\|_G^2$ to denote $\langle x,
Gx \rangle$ if $G$ is a self-adjoint linear operator in a given finite dimensional Euclidean space $\cX$. If $\Sigma:\cX \to \cX$ is
a self-adjoint and positive semi-definite linear operator,  we use
$\Sigma^{\frac{1}{2}}$  to denote the unique self-adjoint and
positive semi-definite  square root of $\Sigma$.

\item Denote
\[
    x: = \left(
          \begin{array}{c}
           x_1        \\
           x_2        \\
           x_3
            \end{array}
           \right),  \qquad u: = \left(
          \begin{array}{c}
           x_2        \\
           x_3
            \end{array}
           \right), \qquad    A := \left(
          \begin{array}{c}
           A_1 \\
           A_2 \\
           A_3
            \end{array}
           \right),  \qquad B := \left(
          \begin{array}{c}
           A_2 \\
           A_3
            \end{array}
           \right).
            \nn \]
\item Let  $\alpha \in (0, 1]$ be given. Denote
\[ \label{MatrixM}
M := \left(
          \begin{array}{cc}
           (1 - \alpha) \Sigma_2 + T_2        &     0   \\
           0                                  &    \Sigma_3 + T_3
            \end{array}
           \right) +   \sigma BB^*\]
\[  \label{MatrixH}
     H := \left(
          \begin{array}{cc}
        \frac{5(1 - \alpha)}{2} \Sigma_2 +  T_2     &     0   \\
           0                                        &   \frac{5}{2}\Sigma_3 + T_3 - \frac{5\sigma^2}{2\alpha}(A_2A_3^*)^*\Sigma_2^{-1}(A_2A_3^*)
            \end{array}
           \right) + \min(\tau, 1+\tau-\tau^2) \sigma BB^*.
             \]
\end{itemize}

\section{A $3$-Block Semi-Proximal ADMM}

Based on our previous introduction and motivation, we propose our  $3$-block semi-proximal
ADMM for solving problem  \eqref{ConvexP-G} in the following:

\bigskip

\centerline{\fbox{\parbox{\textwidth}{ {\bf Algorithm sPADMM: \bf{A
3-block semi-proximal ADMM for solving problem \eqref{ConvexP-G}.}}
\\
Let   $\sigma \in (0, +\infty)$ and $\tau \in (0,
+\infty)$ be given parameters. Let $T_i$, $i=1, 2, 3$,
be given self-adjoint and positive semi-definite linear operators defined on ${\cal
X}_i$, $i=1, 2, 3$, respectively.
 Choose $(x_1^0, x_2^0, x_3^0, z^0)
\in \hbox{dom}(\theta_1) \times \hbox{dom}(\theta_2) \times
\hbox{dom}(\theta_3) \times {\cal Z}$ and set $k=0$.
\begin{description}
\item [Step 1.] Compute
\begin{equation}\label{SPADMM}
\left\{
\begin{array}{l}
 \displaystyle x_1^{k+1} := \argmin_{x_1 \in {\cal X}_1}  \bigl\{
     {\cal L}_{\sigma}(x_1, x_2^k, x_3^k; z^k) + \frac{1}{2}\|x_1 - x_1^k\|_{T_1}^2 \bigr\}, \\
\displaystyle x_2^{k+1} := \argmin_{x_2 \in {\cal X}_2} \bigl\{
     {\cal L}_{\sigma}(x_1^{k+1}, x_2, x_3^k; z^k)+ \frac{1}{2}\|x_2 - x_2^k\|_{T_2}^2 \bigr\}, \\
\displaystyle x_3^{k+1} := \argmin_{x_3 \in {\cal X}_3}  \bigl\{
     {\cal L}_{\sigma}(x_1^{k+1}, x_2^{k+1},  x_3; z^k)
      + \frac{1}{2}\|x_3 - x_3^k\|_{T_3}^2 \bigr\}, \\
z^{k+1} := z^k + \tau\sigma(A^* x^{k+1} - c).
\end{array}
\right.
\end{equation}
\item [Step 2.] If a termination criterion is not met, set $k:=k+1$ and then goto Step 1.
  \end{description}
}}}
\medskip

In order to analyze the convergence properties of Algorithm sPADMM,  we   make the following   assumptions.

\begin{assumption}\label{assump-strong-convex}
The convex  function  $\theta_2$ satisfies \eqref{strong-convex} with $\Sigma_2 \succ 0$.
\end{assumption}
\begin{assumption}\label{assump-strong-convex-B}
The self-adjoint and positive semi-definite operators $T_i$, $i=1,2,3$, are chosen such that the sequence $\{(x_1^k, x_2^k, x_3^k, z^k)\}$ generated by Algorithm  sPADMM is well defined.
\end{assumption}

\begin{assumption}\label{assump-CQ} There exists $x' =(x_1', x_2', x_3') \in
\hbox{ri} (\hbox{dom}(\theta_1)  \times \hbox{dom}(\theta_2) \times
\hbox{dom}(\theta_3)) \bigcap P$, where
\[  P: = \Bigl\{x := (x_1, x_2, x_3) \in {\cal X}_1 \times {\cal X}_2 \times {\cal X}_3 \; \Big| \;   A^*x  = c \Bigr\}. \nn\]
\end{assumption}

Under Assumption \ref{assump-CQ}, it follows from \cite[Corollary
28.2.2]{Rockafellar70} and \cite[Corollary 28.3.1]{Rockafellar70}
that $\bar{x} = (\bar{x}_1, \bar{x}_2, \bar{x}_3) \in {\cal X}_1 \times {\cal
X}_2 \times {\cal X}_3$ is an optimal solution to problem
 \eqref{ConvexP-G} if and only if there exists a Lagrange multiplier $\bar{z} \in
{\cal Z}$ such that
\[\label{gradient-pq}
     -A_i \bar{z} \in  \partial \theta_i(\bar{x}_i), \; i= 1, 2, 3 \quad \hbox{and} \quad
      A^* \bar{x}  - c = 0. \]
 Moreover, any $\bar{z} \in {\cal Z}$
satisfying \eqref{gradient-pq}  is an optimal solution to the dual
 of problem \eqref{ConvexP-G}.

Let $\bar{x} = (\bar{x}_1, \bar{x}_2, \bar{x}_3) \in {\cal X}_1
\times {\cal X}_2 \times {\cal X}_3$  and $\bar z\in \cZ$ satisfy
(\ref{gradient-pq}). For the sake of convenience, define for $(x_1, u,z) :=(x_1,
(x_2, x_3), z) \in {\cal X}_1 \times ({\cal X}_2 \times   {\cal X}_3)
\times {\cal Z}$, $\alpha \in (0, 1]$ and $k = 0, 1, \ldots$, the
following quantities
\begin{equation*}
\begin{array}{lll}
  {\phi}_k (x_1, u, z) &:=&  (\sigma\tau)^{-1} \|z^k - z\|^2
 + \|x_1^k - x_1\|^2_{\Sigma_1 + T_1}
  + \|u^k - u\|^2_M
 \end{array}
\end{equation*}
and
\begin{equation}\label{Notation_dtp1}
\left\{
\begin{array}{l}
x_{ie}^k: =x_i^k - \bar{x}_i, \; i = 1,2,3, \quad  u_{e}^k: =u^k - \bar{u},  \quad z_{e}^k: =z^k - \bar{z},  \\[8pt]
\Delta x_i^k: = x_i^{k+1} - x_i^k, \; i = 1,2,3, \quad  \Delta u^k: = u^{k+1} - u^k,  \quad \Delta z^k: = z^{k+1} - z^k,  \\[8pt]
 \overline{\phi}_k  :=  {\phi}_k(\bar{x}_1, \bar{u}, \bar{z}) = (\sigma\tau)^{-1} \|z_e^k\|^2
 + \|x_{1e}^k \|^2_{\Sigma_1 + T_1}
  + \|u_e^k\|^2_M, \\ [8pt]
\xi_{k+1} : =  \|\Delta x_2^k\|^2_{T_2} + \|\Delta x_3^k\|^2_{T_3 +
\frac{\sigma^2}{\alpha}(A_2A_3^*)^*\Sigma_2^{-1}(A_2A_3^*)}, \\
[8pt]
 s_{k+1} : =   \|\Delta x_1^k\|_{\frac{1}{2}\Sigma_1 + T_1}^2 + \|\Delta
x_2^k\|^2_{\frac{1-\alpha}{2}\Sigma_2 + T_2}  + \|\Delta
x_3^k\|^2_{\frac{1}{2}\Sigma_3 + T_3 -
\frac{\sigma^2}{2\alpha}(A_2A_3^*)^*\Sigma_2^{-1}(A_2A_3^*)} \\[8pt]
\qquad \qquad + \sigma\|A_1^*x_1^{k+1} + B^* u^k - c\|^2, \\[8pt]
 t_{k+1}  : =  \|\Delta x_1^k\|_{\frac{1}{2}\Sigma_1 + T_1}^2 + \|\Delta
u^k\|^2_H, \\[8pt]
   r^k :=  A^* x^k - c.
\end{array}
\right.
\end{equation}

To prove the convergence of Algorithm sPADMM for solving problem \eqref{ConvexP-G}, we
first present some useful  lemmas.

\begin{lemma}  Assume  that  Assumptions \ref{assump-strong-convex}, \ref{assump-strong-convex-B} and \ref{assump-CQ} hold.
Let $\{(x_1^k, x_2^k, x_3^k, z^k)\}$ be generated by Algorithm
sPADMM. Then, for any $\tau \in (0, +\infty)$ and
integer $k \ge 0$, we have
\[\label{com-ine2}
   \overline{\phi}_k  - \overline{\phi}_{k+1}
   \ge   (1 - \tau) \sigma \|r^{k+1}\|^2
   +  s_{k+1},
\]
where $\overline{\phi}_k$, $s_{k+1}$ and $r^{k+1}$ are defined as in \eqref{Notation_dtp1}.
\end{lemma}
\noindent{\bf Proof}. The sequence $\{(x_1^k, x_2^k,
x_3^k, z^k)\}$ is well defined under
Assumption \ref{assump-strong-convex-B}. Notice that the iteration
scheme (\ref{SPADMM})  of Algorithm sPADMM can be re-written as for
$k=0, 1, \ldots$ that
\begin{equation}\label{ADMM-semidefinite}
\left\{
\begin{array}{l}
  - A_1 [z^k + \sigma(A_1^*x_1^{k+1} + \sum_{j=2}^3 A_j^* x_j^k - c)]
  - T_1(x_1^{k+1} - x_1^k) \in \partial \theta_1(x_1^{k+1}), \\[3pt]
  - A_2 [z^k + \sigma(\sum_{j=1}^2 A_j^*x_j^{k+1} +  A_3^* x_3^k - c)]
     - T_2(x_2^{k+1}  - x_2^k) \in \partial \theta_2(x_2^{k+1}), \\[3pt]
  - A_3 [z^k + \sigma( A^* x^{k+1} - c)]
    - T_3(x_3^{k+1}  - x_3^k)\in \partial \theta_3(x_3^{k+1}),
 \\[3pt]
z^{k+1} := z^k + \tau\sigma(A^* x^{k+1} - c).
\end{array}
\right.
\end{equation}
Combining \eqref{strong-convex} with \eqref{gradient-pq} and
\eqref{ADMM-semidefinite}, and using the definitions of
$x_{ie}^{k+1}$ and $\Delta x_i^k$,  for $i=1, 2, 3$, we have
\begin{eqnarray}\label{ADMM-ine-non-1}
 \Big\langle x_{ie}^{k+1}, \;    A_i \bar{z}  - A_i z^k  - \sigma A_i  (\sum_{j=1}^i A_j^* x_j^{k+1} +
\sum_{j=i+1}^3 A_j^* x_j^k - c) - T_i\Delta x_i^k \Big\rangle \ge
\|x_{ie}^{k+1}\|_{\Sigma_i}^2.
\end{eqnarray}
For any vectors $a,b,d$ in the same Euclidean vector  space and any self-adjoint linear  operator $G$, we have the identity
\[
\big\langle a-b, \; G(d-a)\big\rangle = \frac{1}{2}(\|d - b\|^2_G -
\|a - b\|^2_G - \|a - d\|^2_G).
  \nn \]
Taking  $a=x_i^{k+1}$, $b=\bar{x}_i$, $d=x_i^k$  and $G=T_i$ in the above
identity, and using the definitions of $x_{ie}^{k+1}$ and $\Delta
x_i^k$, we get
\[ \label{T-relation}
   \big\langle   x_{ie}^{k+1}, \;  -T_i \Delta x_i^k)\big\rangle
= \frac{1}{2}(\|x_{ie}^k \|^2_{T_i}  - \|x_{ie}^{k+1}\|^2_{T_i} - \|\Delta x_i^k\|^2_{T_i}), \quad i=1,2,3.  \]
Let
\[\label{ADMM-Tnotation}
  {\tilde z}^{k+1}  = z^k + \sigma  (A^* x^{k+1} - c) = z^k + \sigma (A_1^*x_1^{k+1} + B^*u^{k+1} - c). \]
Substituting \eqref{T-relation} and \eqref{ADMM-Tnotation} into
\eqref{ADMM-ine-non-1} and using the definition of $\Delta x_j^k$,
for $i=1, 2$, we have
\[\label{ine-k1}
  \Big\langle   x_{ie}^{k+1}, \; A_i \bar{z}  - A_i \tilde{z}^{k+1} + \sigma A_i \sum_{j=i+1}^3 A_j^* \Delta x_j^k  \Big\rangle
 + \frac{1}{2} (\|x_{ie}^k \|_{T_i}^2 - \|x_{ie}^{k+1}\|_{T_i}^2)  \ge \frac{1}{2} \| \Delta x_i^k\|_{T_i}^2  + \|x_{ie}^{k+1}\|^2_{\Sigma_i}
\]
and
\[\label{ADMM-imp-inem}
   \big\langle    x_{3e}^{k+1}, \; A_3 \bar{z} - A_3 \tilde{z}^{k+1}
  \big\rangle + \frac{1}{2} (\|x_{3e}^k \|^2_{T_3} - \|x_{3e}^{k+1} \|^2_{T_3})
  \ge \frac{1}{2} \|\Delta x_3^k\|^2_{T_3} + \|x_{3e}^{k+1} \|^2_{\Sigma_3}.
\]
Adding \eqref{ine-k1} for $i=1, 2$ to \eqref{ADMM-imp-inem}, we get
\begin{eqnarray}\label{add-equ}
&  &\sum_{i=1}^3 \big\langle   x_{ie}^{k+1}, \; A_i  \bar{z} - A_i  \tilde{z}^{k+1} \big\rangle +
 \sigma  \big\langle    x_{1e}^{k+1}, \;   A_1  \sum_{j=2}^3 A_j^* \Delta x_j^k\big\rangle
 +  \sigma \big\langle    x_{2e}^{k+1}, \;   A_2  A_3^* \Delta x_3^k \big\rangle \nn \\
&  & \;  +   \frac{1}{2}\sum_{i=1}^3 \big(\|x^k_{ie}  \|^2_{T_i} - \|x^{k+1}_{ie} \|^2_{T_i}\big)
   \ge \frac{1}{2} \sum_{i=1}^3\|\Delta x^k_i\|_{T_i}^2 +  \sum_{i=1}^3\|x^{k+1}_{ie} \|_{\Sigma_i}^2.
\end{eqnarray}
By   simple manipulations and using $ A_1^* x_{1e}^{k+1} =
A_1^* x_1^{k+1} - A_1^* \bar{x}_1 = B^*\bar{u} + ( A_1^* x_1^{k+1} -
c)$, we get
\begin{eqnarray}\label{ADMM-0A-1}
  \sigma  \big\langle    x_{1e}^{k+1}, \;   A_1  \sum_{j=2}^3 A_j^* \Delta x_j^k\big\rangle  &  = &
  \sigma  \big\langle  - x_{1e}^{k+1}, \; - A_1 B^* \Delta u^k   \big\rangle
= \sigma  \big\langle   - A_1^* x_{1e}^{k+1}, \;B^* u^k - B^* u^{k+1} \big\rangle \nn \\
&  = & \sigma \big\langle (- B^* \bar{u}) - ( A_1^* x_1^{k+1} - c),
\;
 (- B^* u^{k+1})- (-B^* u^k) \big\rangle.
\end{eqnarray}
For any vectors $a,b,d,e$ in the same Euclidean vector space, we have the identity
\[\label{identity}
\big\langle a-b, \; d-e \big\rangle = \frac{1}{2}(\|a - e\|^2 - \|a
- d\|^2) + \frac{1}{2}(\|b - d\|^2 - \|b - e\|^2).
  \]
In the above  identity, by taking $a=  - B^* \bar{u}$, $b= A_1^* x_1^{k+1} -
c$, $d= - B^* u^{k+1}$ and $e= -B^* u^k$, and applying  it to the
right-hand side of \eqref{ADMM-0A-1}, we obtain  from  the definitions of
$u_e^k$ and $\tilde{z}^{k+1}$ that
\begin{eqnarray}\label{ADMM-0A}
&  &  \sigma  \big\langle    x_{1e}^{k+1}, \;   A_1  \sum_{j=2}^3
A_j^* \Delta x_j^k   \big\rangle \nn \\
&  & \;  =  \frac{\sigma}{2} ( \| B^*
 u_{e}^k  \|^2 -
        \|B^* u_{e}^{k+1}  \|^2 ) +  \frac{\sigma}{2}  (\|A_1^*x_1^{k+1} + B^*u^{k+1} - c\|^2 -
        \| A_1^* x_1^{k+1} + B^* u^k - c \|^2) \nn \\
&  & \; =  \frac{\sigma}{2}  ( \|  B^*
 u_e^k  \|^2 - \|  B^* u_e^{k+1}  \|^2  )  +\frac{1}{2\sigma}  \|z^k - \tilde{z}^{k+1} \|^2 - \frac{\sigma}{2}
        \| A_1^* x_1^{k+1} + B^* u^k - c \|^2.
\end{eqnarray}
Using the Cauchy-Schwarz inequality, for  the parameter $\alpha \in (0, 1]$, we
get
\begin{eqnarray}\label{x2e}
  \sigma \big\langle    x_{2e}^{k+1}, \;   A_2  A_3^* \Delta x_3^k \big\rangle
  & =  & 2\big\langle  (\alpha\Sigma_2)^{\frac{1}{2}}  x_{2e}^{k+1}, \;   \frac{\sigma}{2}(\alpha\Sigma_2)^{-\frac{1}{2}} A_2  A_3^* \Delta x_3^k
  \big\rangle \nn \\
 & \le &  \alpha\|x_{2e}^{k+1}  \|^2_{\Sigma_2} +
\frac{\sigma^2}{4\alpha} \|\Delta x_3^k
\|^2_{(A_2 A_3^*)^*\Sigma_2^{-1}(A_2 A_3^*)}.
 \end{eqnarray}
It follows from \eqref{ADMM-Tnotation} that
\[\label{sum-equ}
  \sum_{i=1}^3 \big\langle   x_{ie}^{k+1}, \; A_i  \bar{z} - A_i  \tilde{z}^{k+1} \big\rangle
   =     \big\langle \bar{z} - \tilde{z}^{k+1}, \; \sum_{i=1}^3 A_i^* x_{ie}^{k+1}  \big\rangle
   =  \frac{1}{\sigma}
   \big\langle \bar{z} - \tilde{z}^{k+1}, \; \tilde{z}^{k+1} - z^k\big\rangle.
\]
Substituting \eqref{ADMM-0A}, \eqref{x2e} and \eqref{sum-equ} into \eqref{add-equ}, we obtain
\begin{eqnarray}\label{com-ine1}
&  &  \frac{1}{\sigma}\big\langle \bar{z} - \tilde{z}^{k+1}, \;
\tilde{z}^{k+1} - z^k\big\rangle + \frac{1}{2\sigma}\|z^k -
\tilde{z}^{k+1}\|^2
  + \frac{\sigma}{2}  ( \|B^* u_e^k \|^2 -  \| B^* u_e^{k+1} \|^2 ) \nn \\
&  & \qquad +   \frac{1}{2}\sum_{i=1}^3 \big(\|x^k_{ie}  \|^2_{T_i} - \|x^{k+1}_{ie} \|^2_{T_i}\big) \nn \\
  &  & \quad \ge \frac{\sigma}{2} \|A_1^* x_1^{k+1} + B^* u^k - c \|^2
  + \frac{1}{2} \sum_{i=1}^3\|\Delta x^k_i\|_{T_i}^2 +  \sum_{i=1, i\neq 2}^3\|x^{k+1}_{ie} \|_{\Sigma_i}^2  \nn \\
  &  & \qquad + (1 - \alpha) \|x^{k+1}_{2e} \|_{\Sigma_2}^2  - \frac{\sigma^2}{4\alpha}
  \| \Delta x_3^k \|^2_{(A_2 A_3^*)^*\Sigma_2^{-1}(A_2 A_3^*)}.
\end{eqnarray}
From the elementary inequality $\|a\|^2 + \|b\|^2 \ge \|a-b\|^2/2$ and $x_{ie}^{k+1} - x_{ie}^k = \Delta x_i^k$,
it follows that
\begin{eqnarray}\label{x123}
 &  &  \sum_{i=1, i\neq 2}^3\|x^{k+1}_{ie} \|_{\Sigma_i}^2 + (1 - \alpha) \|x^{k+1}_{2e}
 \|_{\Sigma_2}^2 \nn \\
 &  & \quad = \frac{1}{2}\sum_{i=1, i\neq 2}^{3}(\|x_{ie}^{k+1} \|^2_{\Sigma_i} + \|x_{ie}^k \|^2_{\Sigma_i})
  + \frac{1}{2}\sum_{i=1, i\neq 2}^{3}(\|x_{ie}^{k+1} \|^2_{\Sigma_i} - \|x_{ie}^k \|^2_{\Sigma_i})  \nn \\
 &  & \qquad + \frac{1-\alpha}{2} (\|x_{2e}^{k+1} \|^2_{\Sigma_2} + \|x_{2e}^k \|^2_{\Sigma_2})
  + \frac{1 - \alpha}{2} (\|x_{2e}^{k+1} \|^2_{\Sigma_2} - \|x_{2e}^k
  \|^2_{\Sigma_2}) \nn \\
&   & \quad \ge  \frac{1}{4}\sum_{i=1, i\neq 2}^{3} \|\Delta
x_i^k\|^2_{\Sigma_i} + \frac{1}{2}\sum_{i=1, i\neq
2}^{3}(\|x_{ie}^{k+1} \|^2_{\Sigma_i}  - \|x_{ie}^k \|^2_{\Sigma_i})
+ \frac{1 - \alpha}{4} \|\Delta x_2^k\|^2_{\Sigma_2} \nn \\
&  & \qquad + \frac{1 - \alpha}{2}(\|x_{2e}^{k+1} \|^2_{\Sigma_2}  -
\|x_{2e}^k \|^2_{\Sigma_2}).
\end{eqnarray}
By   simple manipulations and using the definition of $z_e^k$, we get
\begin{eqnarray}\label{ADMM-lam2}
&  & \frac{1}{\sigma}\big\langle \bar{z} - \tilde{z}^{k+1}, \;
\tilde{z}^{k+1} - z^k \big\rangle+ \frac{1}{2\sigma}\|z^k -
\tilde{z}^{k+1}\|^2 \nn \\
&  &  \quad = \frac{1}{\sigma}\big\langle   \bar{z} - z^k, \; \tilde{z}^{k+1} - z^k\big\rangle
+ \frac{1}{\sigma}\big\langle   z^k - \tilde{z}^{k+1}, \; \tilde{z}^{k+1} - z^k\big\rangle + \frac{1}{2\sigma}\|z^k - \tilde{z}^{k+1}\|^2\nn \\
&  &  \quad = \frac{1}{\sigma}\big\langle   - z_e^k, \; \tilde{z}^{k+1} - z^k\big\rangle - \frac{1}{2\sigma}\|z^k - \tilde{z}^{k+1}\|^2 \nn \\
&  &  \quad = \frac{1}{2\sigma\tau}\Big(\|z_e^k  \|^2 - \|z_e^k   +
\tau(\tilde{z}^{k+1} - z^k)\|^2\Big) +
             \frac{\tau - 1}{2\sigma} \|z^k - \tilde{z}^{k+1}\|^2.
\end{eqnarray}
By using \eqref{ADMM-semidefinite}, \eqref{ADMM-Tnotation} and the definitions of $z_e^k$ and $r^{k+1}$, we have
\[
 z_e^{k+1} = z_e^k +   \tau (\tilde{z}^{k+1} - z^k) \qquad
 \hbox{and} \qquad z^k -  \tilde{z}^{k+1} = - \sigma r^{k+1}, \nn \]
which, together with \eqref{ADMM-lam2}, imply
\[ \label{ADMM-lam3} \frac{1}{\sigma}\big\langle \bar{z} - \tilde{z}^{k+1}, \; \tilde{z}^{k+1} - z^k\big\rangle + \frac{1}{2\sigma}\|z^k - \tilde{z}^{k+1}\|^2
 = \frac{1}{2\sigma\tau}(\|z_e^k \|^2 - \|z_e^{k+1} \|^2)
 + \frac{(\tau - 1)\sigma}{2}\|r^{k+1}\|^2. \]
Substituting \eqref{x123} and \eqref{ADMM-lam3} into
\eqref{com-ine1}, and using the definitions of $\overline{\phi}_k$,
$s_{k+1}$ and $r^{k+1}$, we get the assertion \eqref{com-ine2}. The
proof is complete.   \hfill {$\Box$}

\begin{lemma}\label{lem-ine-new}  Assume  that  Assumptions \ref{assump-strong-convex} and \ref{assump-strong-convex-B} hold.
 Let $\{(x_1^k, x_2^k, x_3^k, z^k)\}$ be generated by Algorithm sPADMM. Then,
for any $\tau \in (0, +\infty)$ and integer $k \ge 1$, we have
\begin{eqnarray}\label{ADMM-ine-10}
&  & - \sigma \big\langle B^* \Delta u^k, \;  r^{k+1}  \big\rangle
 \ge - (1 - \tau)\sigma  \big\langle B^*  \Delta u^k, \; r^k  \big\rangle +
\frac{1}{2}\sum_{i=2}^3 ( \|\Delta x_i^k\|^2_{T_i + 2\Sigma_i} -
\|\Delta x_i^{k-1} \|^2_{T_i})
              \nn \\
&  &  \qquad + \sigma  \big\langle A_2^* \Delta x_2^k, \; A_3^*
(\Delta x_3^{k-1}    -
  \Delta x_3^k) \big\rangle,
\end{eqnarray}
where $\Delta u^k$, $\Delta x_i^k$ $(i=2,3)$ and $r^{k+1}$ are defined as in \eqref{Notation_dtp1}.
\end{lemma}
\noindent{\bf Proof}. Let
\[ v^{k+1}: = z^k +
\sigma\Big(\sum_{j=1}^2 A_j^* x_j^{k+1} +
    A_3^* x_3^k - c\Big). \nn \]
By using \eqref{ADMM-semidefinite} and   the definition
of $\Delta x_2^k$, we have
\[  - A_2  v^{k+1} - T_2 \Delta x_2^k  \in \partial \theta_2(x_2^{k+1}) \quad
   \hbox{and} \quad - A_2  v^k - T_2\Delta x_2^{k-1}  \in \partial \theta_2(x_2^k).\nn\]
Thus, we obtain from \eqref{strong-convex} that
\[
 \big\langle   \Delta x_2^k, \; (A_2 v^k + T_2\Delta x_2^{k-1}) - (A_2 v^{k+1} + T_2 \Delta x_2^k)
 \big\rangle\ge \|\Delta x_2^k\|^2_{\Sigma_2}. \nn
\]
By using the Cauchy-Schwarz inequality, we obtain
\[
  \big\langle \Delta x_2^k, \; T_2(\Delta x_2^k  - \Delta x_2^{k-1})
 \big\rangle
 =   \|\Delta x_2^k\|^2_{T_2} - \big\langle \Delta x_2^k, \;  T_2\Delta x_2^{k-1}
 \big\rangle   \ge  \frac{1}{2} \|\Delta x_2^k\|^2_{T_2}    - \frac{1}{2} \|\Delta x_2^{k-1}\|^2_{T_2}. \nn
\]
Adding up the above two inequalities, we get
\[ \label{ADMM-ine-2}
  \big\langle  A_2^* \Delta x_2^k, \; v^k - v^{k+1}\big\rangle
   \ge \frac{1}{2}\|\Delta x_2^k\|_{T_2 + 2\Sigma_2}^2 - \frac{1}{2}\|\Delta x_2^{k-1}\|_{T_2}^2. \]
Using $z^{k-1} - z^k = - \tau\sigma r^k$ and the definitions of
$v^k$ and $r^k$, we have
\[   v^k - v^{k+1} = (1 - \tau)\sigma r^k - \sigma r^{k+1} - \sigma A_3^* (\Delta x_3^{k-1} - \Delta x_3^k). \nn\]
Substituting the above equation  into \eqref{ADMM-ine-2}, we get
\begin{eqnarray}\label{ADMM-ine-3-1}
  \sigma \big\langle - A_2^* \Delta x_2^k, \;   r^{k+1} \big\rangle    & \ge & - (1 - \tau)\sigma \big\langle A_2^* \Delta x_2^k, \;
r^k  \big\rangle
              + \sigma  \big\langle A_2^* \Delta x_2^k, \;   A_3^*  (\Delta x_3^{k-1} - \Delta x_3^k)\big\rangle  \nn \\
&  &
 + \frac{1}{2}  (\|\Delta x_2^k\|^2_{T_2 + 2\Sigma_2} -  \|\Delta x_2^{k-1}\|^2_{T_2}).
\end{eqnarray}
Similarly as for deriving \eqref{ADMM-ine-3-1}, we can obtain that
\[
 \sigma  \big\langle - A_3^* \Delta x_3^k,  r^{k+1}  \big\rangle  \ge - (1 - \tau)\sigma \big\langle A_3^* \Delta x_3^k,
r^k  \big\rangle
 + \frac{1}{2} ( \|\Delta x_3^k\|^2_{T_3 + 2\Sigma_3} -   \|\Delta x_3^{k-1}\|^2_{T_3}).
\nn\]
Adding up the above
inequality and   \eqref{ADMM-ine-3-1},   and using the definitions of $B^*$ and $u$, we get the
assertion \eqref{ADMM-ine-10}. The proof is complete. \hfill
{$\Box$}

\begin{lemma}\label{lem-ine}  Assume  that  Assumptions \ref{assump-strong-convex} and \ref{assump-strong-convex-B} hold.
Let $\{(x_1^k, x_2^k, x_3^k, z^k)\}$ be generated by Algorithm
sPADMM. For any $\tau \in (0, +\infty)$ and integer
$k \ge 1$, we have
\begin{eqnarray}\label{mul-ine-0}
  (1 - \tau) \sigma\|r^{k+1}\|^2 +  s_{k+1}
   & \ge &  t_{k+1}  +  \max(1 - \tau, 1 - \tau^{-1})\sigma(\|r^{k+1}\|^2 - \|r^k\|^2)
  \nn \\
 &  &  +  \min(\tau,  1+ \tau - \tau^2)\sigma \tau^{-1} \|r^{k+1}\|^2 +  (\xi_{k+1} - \xi_k),
\end{eqnarray}
where $s_{k+1}$, $t_{k+1}$, $\xi_{k+1}$ and $r^{k+1}$ are defined as in \eqref{Notation_dtp1}.
\end{lemma}
\noindent{\bf Proof}. By  simple manipulations and using the
definition of $r^{k+1}$, we obtain
\[\label{div-ine}
    \|A_1^* x_1^{k+1} + B^* u^k - c \|^2   =   \|r^{k+1}  -  B^* \Delta u^k  \|^2   =     \|  r^{k+1}  \|^2 - 2
 \big\langle B^* \Delta u^k, \; r^{k+1} \big \rangle  +   \| B^* \Delta
u^k \|^2.
\]
It follows from \eqref{ADMM-ine-10} and \eqref{div-ine} that
\begin{eqnarray}\label{ADMM-ine6}
&  &  (1 - \tau) \sigma\|r^{k+1}\|^2 +
 \sigma \|A_1^* x_1^{k+1} + B^* u^k -
c \|^2 \nn \\
& &   \ge \sigma  \|B^* \Delta u^k \|^2 +  (2-\tau)\sigma  \|r^{k+1}
 \|^2  - 2(1 - \tau)\sigma  \big\langle B^* \Delta u^k, \; r^k  \big\rangle
              + 2\sigma  \big\langle A_2^* \Delta x_2^k, \; A_3^* (\Delta x_3^{k-1} - \Delta x_3^k) \big\rangle \nn \\
&  &  \quad   +
 \sum_{i=2}^3 \big( \|\Delta x_i^k\|^2_{T_i+2\Sigma_i} -  \|\Delta x_i^{k-1} \|^2_{T_i}\big).
\end{eqnarray}
By the Cauchy-Schwarz inequality, for the parameter $\alpha \in (0, 1]$, we
have
\begin{eqnarray*}
&  & 2\sigma  \big\langle A_2^* \Delta x_2^k, \; A_3^* (\Delta x_3^{k-1} - \Delta x_3^k) \big\rangle   \nn \\
&  & \quad = 2\big\langle (\alpha\Sigma_2)^{\frac{1}{2}}\Delta
x_2^k, \; \sigma(\alpha\Sigma_2)^{-\frac{1}{2}} (A_2 A_3^*)  \Delta
x_3^{k-1} \big\rangle
   - 2\big\langle (\alpha\Sigma_2)^{\frac{1}{2}}\Delta x_2^k, \; \sigma(\alpha\Sigma_2)^{-\frac{1}{2}} (A_2 A_3^*)  \Delta x_3^k \big\rangle\nn \\
&  & \quad \ge - \ \alpha\| \Delta x_2^k  \|^2_{\Sigma_2} -
\frac{\sigma^2}{\alpha} \| \Delta x_3^{k-1}  \|^2_{(A_2
A_3^*)^*\Sigma_2^{-1}(A_2 A_3^*)} - \alpha\|\Delta x_2^k
 \|^2_{\Sigma_2}   -    \frac{\sigma^2}{\alpha} \| \Delta x_3^k  \|^2_{(A_2 A_3^*)^*\Sigma_2^{-1}(A_2 A_3^*)}   \nn \\
&  & \quad = - 2\alpha \| \Delta x_2^k \|^2_{\Sigma_2} -
\frac{\sigma^2}{\alpha} (\| \Delta x_3^{k-1} \|^2_{(A_2
A_3^*)^*\Sigma_2^{-1}(A_2 A_3^*)} + \| \Delta x_3^k  \|^2_{(A_2
A_3^*)^*\Sigma_2^{-1}(A_2 A_3^*)}).
\end{eqnarray*}
Substituting the above inequality into \eqref{ADMM-ine6}, we get
\begin{eqnarray}\label{ADMM-ine-new}
&  &  (1 - \tau) \sigma\|r^{k+1}\|^2 +
 \sigma \|A_1^* x_1^{k+1} + B^* u^k -
c \|^2 \nn \\
& &   \ge \sigma  \|B^* \Delta u^k \|^2 +  (2-\tau)\sigma  \|r^{k+1}
 \|^2  - 2(1 - \tau)\sigma  \big\langle B^* \Delta u^k, \; r^k  \big\rangle
  +  \big( \|\Delta x_2^k\|^2_{T_2} -  \|\Delta x_2^{k-1} \|^2_{T_2}\big)             \nn \\
&  &  \quad   +
  \big( \|\Delta x_3^k\|^2_{T_3+\frac{\sigma^2}{\alpha}(A_2 A_3^*)^*\Sigma_2^{-1}(A_2 A_3^*)}
  -  \|\Delta x_3^{k-1} \|^2_{T_3+\frac{\sigma^2}{\alpha}(A_2 A_3^*)^*\Sigma_2^{-1}(A_2 A_3^*)}\big)
   + 2(1-\alpha)\|\Delta x_2^k\|^2_{\Sigma_2} \nn \\
&  & \quad + 2 \|\Delta x_3^k\|^2_{\Sigma_3} -
\frac{2\sigma^2}{\alpha}  \|\Delta x_3^k\|^2_{(A_2
A_3^*)^*\Sigma_2^{-1}(A_2 A_3^*)}.
\end{eqnarray}
By using the definitions of $s_{k+1}$ and $t_{k+1}$, and the fact
that
\[ \|\Delta u^k\|^2_H =   \| \Delta x_2^k  \|^2_{\frac{5(1-\alpha)}{2}\Sigma_2 + T_2}
+ \|\Delta x_3^k\|^2_{\frac{5}{2}\Sigma_3 + T_3 - \frac{5\sigma^2}{2\alpha}(A_2A_3^*)^*\Sigma_2^{-1}(A_2A_3^*)}
+ \min(\tau,  1+ \tau - \tau^2)\sigma \|B^* \Delta u^k\|^2,\nn\] we
have
\begin{eqnarray*}
&  &  2(1-\alpha)\|\Delta x_2^k\|^2_{\Sigma_2}  + 2 \|\Delta x_3^k\|^2_{\Sigma_3}
- \frac{2\sigma^2}{\alpha}  \|\Delta x_3^k\|^2_{(A_2 A_3^*)^*\Sigma_2^{-1}(A_2 A_3^*)} \nn \\
&  & \quad =   - s_{k+1} + t_{k+1} -  \min(\tau,  1+ \tau -
\tau^2)\sigma \|B^* \Delta u^k\|^2 + \sigma \|A_1^* x_1^{k+1} + B^*
u^k - c \|^2.
\end{eqnarray*}
Substituting the above equation into \eqref{ADMM-ine-new} and using
the definition of $\xi_{k+1}$, we get
\begin{eqnarray}\label{ADMM-ine-new2}
&  &  (1 - \tau) \sigma\|r^{k+1}\|^2 +   s_{k+1} -
t_{k+1} +    \min(\tau,  1+ \tau - \tau^2)\sigma \|B^* \Delta u^k\|^2 \nn \\
&  & \quad  \ge  \sigma \|B^*  \Delta u^k \|^2 +  (2-\tau)\sigma \|r^{k+1}  \|^2
     -   2(1 - \tau)\sigma  \big\langle B^*
 \Delta u^k, \; r^k \big\rangle +  (\xi_{k+1} -
\xi_k).
\end{eqnarray}
By using the Cauchy-Schwarz inequality, we get
\begin{equation}\label{tau-Cauchy}
\left\{
\begin{array}{ll}
 - 2(1 - \tau)\sigma  \big\langle B^* \Delta u^k, \;r^k  \big\rangle \ge - (1 - \tau)\sigma \|B^*
\Delta u^k\|^2 -  (1 - \tau)\sigma \|r^k\|^2 & \hbox{if} \; \tau \in (0,1], \\[8pt]
- 2(1 - \tau)\sigma  \big\langle B^* \Delta u^k, \;   r^k
\big\rangle \ge  (1 - \tau)\tau\sigma \|B^* \Delta u^k\|^2 +
\frac{(1 - \tau)\sigma}{\tau}\| r^k \|^2 &  \hbox{if} \; \tau \in
(1, +\infty).
\end{array}
\right.
\end{equation}
Substituting \eqref{tau-Cauchy} into \eqref{ADMM-ine-new2},  we obtain from
simple manipulations that
\begin{eqnarray*}
&  &  (1 - \tau) \sigma\|r^{k+1}\|^2 +   s_{k+1} -
t_{k+1} +    \min(\tau,  1+ \tau - \tau^2)\sigma \|B^* \Delta u^k\|^2 \nn \\
&  & \quad  \ge  \max(1 - \tau, 1 - \tau^{-1})\sigma(\|r^{k+1}\|^2 - \|r^k\|^2)
+  \min(\tau,  1+ \tau - \tau^2)\sigma (\tau^{-1} \|r^{k+1}\|^2 + \|B^* \Delta u^k\|^2)\\
&  & \qquad +  (\xi_{k+1} - \xi_k).
\end{eqnarray*}
The assertion \eqref{mul-ine-0} is proved immediately.
 \hfill {$\Box$}

Now, we are ready to prove the   convergence of the sequence
$\{(x_1^k, x_2^k, x_3^k, z^k)\}$ generated by Algorithm sPADMM.

\begin{theorem}\label{Conver-Alg}
Assume that
 Assumptions \ref{assump-strong-convex}, \ref{assump-strong-convex-B} and \ref{assump-CQ} hold.   Let $\{(x_1^k, x_2^k, x_3^k,
z^k)\}$ be generated by Algorithm sPADMM. Then,  for any
$\tau \in (0, +\infty)$ and integer $k \ge 1$, we have
\begin{eqnarray}\label{case-II2}
&  &   \big(\overline{\phi}_k + \max(1-\tau, 1 - \tau^{-1})\sigma
\|r^k\|^2 + \xi_k\big)  - \big(\overline{\phi}_{k+1}  +
\max(1-\tau, 1 -
\tau^{-1})\sigma  \|r^{k+1}\|^2 + \xi_{k+1}\big) \nn\\
&  &   \quad \ge    t_{k+1} +   \min(\tau,  1+ \tau -
\tau^2)\sigma\tau^{-1} \|r^{k+1}\|^2,
\end{eqnarray}
where $\overline{\phi}_k$, $\xi_{k+1}$, $t_{k+1}$ and $r^k$ are
defined as in \eqref{Notation_dtp1}. Assume that  $\tau \in (0,
(1 + \sqrt{5})/2)$.    If  for some $\alpha \in (0,1]$ it holds that
\begin{equation}\label{eq:sufficient_condition}
\frac{1}{2}\Sigma_1 + T_1 + \sigma A_1 A_1^* \succ
0, \quad H \succ 0 \quad {\rm and} \quad M \succ 0,
\end{equation}
then   the whole sequence $\{(x_1^k, x_2^k, x_3^k)\}$
converges to an optimal solution to problem \eqref{ConvexP-G} and $\{z^k\}$
converges to an optimal solution to the dual   of problem
\eqref{ConvexP-G}.
\end{theorem}
\noindent{\bf Proof}. \noindent  By  substituting \eqref{mul-ine-0} into
\eqref{com-ine2}, we can easily get \eqref{case-II2}.

Assume  that $\tau \in (0, (1 + \sqrt{5})/2)$. Since
\eqref{eq:sufficient_condition} holds for some $\alpha \in (0,1]$, we have $\min(\tau,
1+\tau-\tau^2) >0$, $H \succ 0$ and $M \succ 0$. From
\eqref{case-II2}, we see immediately that the sequence
$\{\overline{\phi}_{k+1}\}$ is bounded, $\lim_{k \rightarrow \infty}
t_{k+1} = 0$ and $\lim_{k \rightarrow \infty}\|r^{k+1}\|=0$, i.e.,
\[\label{limt}
 \lim_{k \rightarrow \infty} \|\Delta x_1^k\|^2_{\frac{1}{2}{\Sigma}_1 + T_1} = 0, \quad
    \lim_{k \rightarrow \infty} \|\Delta u^k\|^2_H = 0, \quad
  \lim_{k \rightarrow \infty}\|r^{k+1}\| = \lim_{k \rightarrow \infty}(\tau\sigma)^{-1}\|\Delta z^k\|
    = 0.
\]
Since $H \succ 0$, we also have that
\[ \label{byy}
   \lim_{k \rightarrow \infty} \|\Delta x_2^k \| = 0,
 \qquad \quad \lim_{k \rightarrow \infty} \|\Delta x_3^k \| = 0
\]
and thus
\[ \label{byy-1}\|A_1^* \Delta x_1^k \| = \Big\|r^{k+1} - r^k - (\sum_{j=2}^3A_j^*\Delta x_j^k) \Big\|\; \le\;  \|r^{k+1}\| + \|r^k\|
     + \sum_{j=2}^3\|A_j^*\Delta x_j^k \| \rightarrow 0 \]
as $k \rightarrow \infty$.
Now from \eqref{limt} and \eqref{byy-1}, we obtain
\[\label{limxk}
  \lim_{k \rightarrow \infty} \|\Delta x_1^k\|_{(\frac{1}{2}\Sigma_1  + T_1 + \sigma A_1 A_1^*)}^2
 = \lim_{k \rightarrow \infty}\big( \|\Delta x_1^k\|_{\frac{1}{2}\Sigma_1  + T_1}^2
       + \sigma \|A_1^*\Delta x_1^k \|^2 \big) = 0.
\]
Recall that $\frac{1}{2}\Sigma_1  + T_1 + \sigma A_1  A_1^* \succ
0$. Thus it follows from \eqref{limxk} that
\[ \label{x-lim}  \lim_{k \rightarrow \infty}\|\Delta x_1^k\| = 0.
\]
By the definition of $\overline{\phi}_{k+1}$, we see that the
three sequences $\{\|z_e^{k+1} \|\}$,  $\{\|x_{1e}^{k+1} \|_{\Sigma_1 + T_1}\}$, and $\{\|u_e^{k+1} \|_M\}$ are all
bounded. Since $M \succ 0$, the sequences $\{\|x_2^{k+1}\|\}$ and $\{\|x_3^{k+1}\|\}$ are also bounded.
Furthermore, by using
\[\label{Axdiff}
\|A_1^* x_{1e}^{k+1}  \|  = \Big\| A^* x^{k+1}  -
 A^* \bar{x} - B^* u_e^{k+1}  \Big\|
\le  \|r^{k+1}\| + \|B^* u_e^{k+1} \|,
\]
we also know that the sequence $\{\|A_1^* x_{1e}^{k+1} \|\}$ is
bounded, and so is the sequence $\{\|x_{1e}^{k+1} \|_{(\Sigma_1 +
T_1 + \sigma A_1 A_1^*)}\}$. This shows that the sequence
$\{\|x_1^{k+1}\|\}$ is also bounded as the operator $\Sigma_1  + T_1
+ \sigma A_1 A_1^* \succeq \frac{1}{2}\Sigma_1  + T_1 + \sigma A_1
A_1^* \succ 0.$
 Thus, the sequence $\{(x_1^k, x_2^k, x_3^k, z^k)\}$ is bounded.

Since the sequence $\{(x_1^k, x_2^k, x_3^k, z^k)\}$ is bounded, there is a subsequence $\{(x_1^{k_i}, x_2^{k_i}, x_3^{k_i}, z^{k_i})\}$
which converges to a cluster point, say $\{(x_1^{\infty}, x_2^{\infty}, x_3^{\infty}, z^{\infty})\}$.
Taking limits on both sides of \eqref{ADMM-semidefinite} along
the subsequence $\{(x_1^{k_i}, x_2^{k_i}, x_3^{k_i}, z^{k_i})\}$, using
\eqref{limt}, \eqref{byy} and \eqref{x-lim},  we obtain that
\begin{equation*}
 - A_j z^{\infty} \in  \partial {\theta}_j(x_j^{\infty}), \; j = 1, 2, 3 \quad \hbox{and} \quad
      A^* x^{\infty} - c = 0,
\end{equation*}
i.e., $(x_1^{\infty}, x_2^{\infty}, x_3^{\infty}, z^{\infty})$ satisfies \eqref{gradient-pq}.
Thus $\{(x_1^{\infty}, x_2^{\infty}, x_3^{\infty})\}$ is an optimal solution
to \eqref{ConvexP-G} and ${z^{\infty}}$ is an optimal solution to the dual of  problem  \eqref{ConvexP-G}.

To complete the proof, we show next that $(x_1^{\infty},
x_2^{\infty}, x_3^{\infty}, z^{\infty})$  is actually the unique
limit of $\{(x_1^k, x_2^k, x_3^k, z^k)\}$. Replacing  $(\bar x_1,
\bar u, \bar z): = (\bar x_1, (\bar x_2, \bar x_3), \bar z)$ by
$(x_1^{\infty}, u^{\infty},  z^{\infty}): =(x_1^{\infty},
(x_2^{\infty}, x_3^{\infty}), z^{\infty})$ in \eqref{case-II2}, for
any integer $k \ge k_i$, we have
\begin{eqnarray}
\label{theta-bound2}
 &  & {\phi}_{k+1}(x_1^{\infty}, u^{\infty}, z^{\infty})  + \max(1-\tau, 1 - \tau^{-1})\sigma
\|r^{k+1}\|^2 + \xi_{k+1}  \nn \\
&  & \quad \le
   {\phi}_{k_i}(x_1^{\infty}, u^{\infty}, z^{\infty}) + \max(1-\tau, 1 - \tau^{-1})\sigma
\|r^{k_i}\|^2 + \xi_{k_i}.
\end{eqnarray}
Note that
$$ \lim_{i \rightarrow \infty} \big({\phi}_{k_i}(x_1^{\infty}, u^{\infty},  z^{\infty}) + \max(1-\tau, 1 - \tau^{-1})\sigma
\|r^{k_i}\|^2 + \xi_{k_i}\big) = 0. $$
 Therefore, from \eqref{theta-bound2} we get
\[ \lim_{k \rightarrow \infty} {\phi}_{k+1}(x_1^{\infty}, u^{\infty}, z^{\infty}) = 0,  \nn \]
i.e.,
\[ \lim_{k \rightarrow \infty} \big( (\sigma\tau)^{-1}\|z^{k+1} - z^{\infty}\|^2 + \|x_1^{k+1} - x_1^{\infty}\|^2_{\Sigma_1+T_1} +
         \|u^{k+1} - u^{\infty}\|^2_M\big)= 0.
  \nn \]
Since $M \succ 0$, we also have that
$\lim_{k \rightarrow \infty} u^k = u^{\infty}$, that is $\lim_{k \rightarrow \infty} x_2^k = x_2^{\infty}$
and $\lim_{k \rightarrow \infty} x_3^k = x_3^{\infty}$. Using the fact that
$\lim_{k\rightarrow\infty} \|r^{k+1}\| = 0$ and
$\lim_{k\rightarrow\infty} \|u^{k+1}-u^{\infty}\| = 0$, we get from
\eqref{Axdiff} that $\lim_{k\rightarrow\infty} \|A_1^*(x_1^{k+1}-x_1^{\infty})\| = 0$.
Thus
\begin{eqnarray*}
& \lim_{k \rightarrow \infty}  \|x_1^{k+1} - x_1^{\infty}\|^2_{\Sigma_1 + T_1  +
     \sigma A_1 A_1^*} = 0. &
\end{eqnarray*}
Since $\Sigma_1 + T_1  + \sigma A_1 A_1^* \succ 0$, we also obtain
that $\lim_{k \rightarrow \infty} x_1^k = x_1^{\infty}$. Therefore,
we have shown that  the sequence $\{(x_1^k, x_2^k, x_3^k)\}$
converges to an optimal solution to \eqref{ConvexP-G} and $\{z^k\}$
converges to an optimal solution to the dual of  problem
\eqref{ConvexP-G} for any $\tau \in (0, (1 + \sqrt{5})/2)$. The proof is complete. \hfill
$\Box$

\begin{remark}\label{condition-example}
Assume that $(1 - \alpha) \Sigma_2 + \sigma A_2A_2^*$ is invertible
for some $\alpha \in (0,1]$. Set $\tau = 1$ (the case that $1\ne
\tau \in (0, (1+\sqrt{5})/2)$ can be discussed in a similar but
slightly more complicated manner)  and $T_2 = 0$ in \eqref{MatrixM}
and \eqref{MatrixH}. Then the assumptions $H \succ 0$ and $M \succ
0$ in \eqref{eq:sufficient_condition} reduce  to
\[   \left(
          \begin{array}{cc}
\frac{5(1 - \alpha)}{2} \Sigma_2 + \sigma A_2A_2^*  &     \sigma A_2A_3^*   \\
           \sigma A_3A_2^*                          &   \frac{5}{2}\Sigma_3 + T_3  +  \sigma A_3A_3^* - \frac{5\sigma^2}{2\alpha}(A_2A_3^*)^*\Sigma_2^{-1}(A_2A_3^*)
            \end{array}  \right) \succ 0  \nn\]
and
\[ \left(
          \begin{array}{cc}
 (1 - \alpha) \Sigma_2 + \sigma A_2A_2^*   &     \sigma A_2A_3^*   \\
           \sigma A_3A_2^*                 &    \Sigma_3 + T_3  +  \sigma A_3A_3^*
            \end{array}  \right) \succ 0, \nn\]
which are, respectively,  equivalent to
\[ \label{cond-1} \frac{5}{2}\Sigma_3 + T_3  +  \sigma A_3A_3^*- \frac{5\sigma^2}{2\alpha}(A_2A_3^*)^*\Sigma_2^{-1}(A_2A_3^*)
    - \sigma^2(A_3A_2^*)\Big(\frac{5(1 - \alpha)}{2} \Sigma_2 + \sigma A_2A_2^*\Big)^{-1}(A_2A_3^*) \succ 0  \]
and
\[ \label{cond-2} \Sigma_3 + T_3   +  \sigma A_3A_3^*
    - \sigma^2(A_3A_2^*)\big((1 - \alpha) \Sigma_2 + \sigma A_2A_2^*\big)^{-1}(A_2A_3^*) \succ 0 \]
    in terms of the Schur-complement format.
The conditions \eqref{cond-1} and \eqref{cond-2} can be satisfied
easily by choosing a proper $T_3$ for  given $\alpha \in (0,1]$ and
$\sigma \in (0, +\infty)$. Evidently, with a fixed $\alpha$, $T_3$
can take a smaller value with  a smaller $\sigma$ and $T_3$ can even
take the zero operator for any $\sigma>0$ smaller than a certain
threshold if  $   \Sigma_3 +   (1-\alpha)\sigma
   A_3A_3^*\succ 0$.  To see this, let us consider the
following  example  constructed  in \cite{ChenHeYeYuan}:
\begin{eqnarray}\label{example-2}
  \begin{array}{rllll}
    \min &   \displaystyle  \frac{1}{20}x_1^2 + \frac{1}{20}x_2^2 + \frac{1}{20}x_3^2    \\[8pt]
   \mbox{\rm s.t.} &  \left(
          \begin{array}{ccc}
           1     &     1     &   1   \\
           1     &     1     &   2   \\
           1     &     2     &   2
            \end{array}
           \right) \left(
          \begin{array}{c}
           x_1   \\
           x_2   \\
           x_3
            \end{array}
           \right) = 0,
   \end{array}
\end{eqnarray}
which is a convex minimization problem with three strongly convex
functions. In \cite{ChenHeYeYuan}, Chen, He, Ye and Yuan showed that
the directly  extended $3$-block ADMM scheme \eqref{extend-ADMM}
with $\tau=\sigma=1$ applied to problem \eqref{example-2}  is
divergent. For problem \eqref{example-2}, $\Sigma_1 = \Sigma_2=
\Sigma_3 = \frac{1}{10}$, $A_1 = (1, 1, 1)$, $A_2 = (1, 1, 2)$ and
$A_3 = (1, 2, 2)$. From \eqref{cond-1} and \eqref{cond-2}, by taking
$\alpha=1$, we have that $T_3$ and $\sigma$ should satisfy the
following conditions
\[ \frac{1}{4} + T_3 - 1225 \sigma^2 + \frac{5}{6}\sigma > 0
  \quad \hbox{and} \quad
  \frac{1}{10} + T_3   + \frac{5}{6}\sigma > 0,  \nn \]
which hold true, in particular, if  $T_3 = 0$ and $\sigma < \frac{1+\sqrt{1765}}{2940} \approx 0.015$ or if $\sigma=1$
and  $T_3 > \frac{14687}{12} \approx 1223.92$.
\end{remark}

\begin{remark} \label{condition-comment} If $A_2^*$ is vacuous, then for any integer $k \ge 0$, we have that
$x_2^{k+1} = x_2^0 = \bar{x}_2$, the $3$-block sPADMM
is just a $2$-block sPADMM, and condition
\eqref{eq:sufficient_condition} reduces  to
\[  \frac{1}{2}\Sigma_1 + T_1 + \sigma A_1  A_1^* \succ 0, \quad
 \Sigma_3 + T_3 + \sigma A_3 A_3^*\succ
0 \quad \hbox{and} \quad \frac{5}{2}\Sigma_3 + T_3 + \min(\tau, 1+ \tau -
\tau^2)\sigma A_3 A_3^*\succ 0, \nn \]
which is equivalent to
\[ \label{2block-condition} \Sigma_1 + T_1 + \sigma A_1  A_1^* \succ 0  \quad
  \hbox{and} \quad \Sigma_3 + T_3 +  \sigma A_3 A_3^*\succ 0   \]
since $\Sigma_1 \succeq 0$, $T_1 \succeq 0$, $\Sigma_3 \succeq 0$ and $T_3 \succeq 0$.
Condition  \eqref{2block-condition} is exactly  the same as the one used in Theorem B.1.
in \cite{FsPoSunTse}.
\end{remark}

\section{Conclusions}

In this paper, we  provided a convergence analysis about  a $3$-block
semi-proximal ADMM for solving  separable convex minimization
problems with the condition  that  the second   block in the objective is
  strongly convex. The step-length $\tau$ in our
proposed semi-proximal ADMM  is allowed to stay in the desirable
region $(0, (1 + \sqrt{5})/2)$. From Remark \ref{condition-example},
we know that with a fixed parameter $\alpha \in (0,1]$, the added
semi-proximal terms can be chosen to be small if the penalty
parameter $\sigma$ is small. If  $A_1^*$ and $A_3^*$ are both
injective and $\sigma>0$ is taken to be smaller than a certain
threshold,  then  the convergent $3$-block semi-proximal ADMM
includes  the directly extended $3$-block ADMM with  $\tau \in (0,
(1+\sqrt{5})/2)$ by taking $T_i$, $i=1,2,3$, to be zero operators.
With no much difficulty, one could extend our $3$-block
semi-proximal ADMM to deal with  the $m$-block  ($m\ge 4$) separable
convex minimization problems possessing $m-2$ strongly convex blocks
and provide the iteration complexity analysis for the corresponding
algorithm in the sense of  \cite{HeYuanSIAM12}. In this work,  we
choose not to do the extension because we are not aware of
interesting applications of  the  $m$-block  ($m\ge 4$) separable
convex minimization problems with $m-2$ strongly convex blocks.
While our sufficient condition bounding the range of values for
$\sigma$ and $T_3$ is quite flexible, it may have one potential
limitation: $T_3$ can be very large if $\sigma$ is not small as
shown  in Remark \ref{condition-example}.  Since a larger $T_3$ can
potentially make the algorithm converge slower, we do not feel that
the study on the iteration complexity is of significance at the
moment unless of course the above  potential limitation is
completely circumvented.

\end{document}